


\input amstex

\documentstyle{amsppt}

\loadbold

\magnification=\magstep1
\parindent=1em
\baselineskip 15pt
\hsize=12.3cm
\vsize=18.5cm

\def\max{\operatorname{max}}
\def\prim{\operatorname{prim}}
\def\ann{\operatorname{ann}}
\def\Rspace{\text{$R$-space}}
\def\Aspace{\text{$A$-space}}
\def\Spec{\operatorname{Spec}}
\def\Mod{\operatorname{Mod}}
\def\Gab{1}
\def\GooWar{2}
\def\Kra{3}
\def\KraSao{4}
\def\McCRobi{5}
\def\McCRobii{6}
\def\Ros{7}
\def\Smi{8}
\def\SmiZan{9}
\def\StaVdB{10}
\def\Ste{11}
\def\VdB{12}

\topmatter

\title Noetherianity of the Space of Irreducible Representations
\endtitle

\rightheadtext{The Space of Irreducible Representations}

\author Edward S. Letzter \endauthor

\abstract Let $R$ be an associative ring with identity. We study an elementary
generalization of the classical Zariski topology, applied to the set
of isomorphism classes of simple left $R$-modules (or, more generally,
simple objects in a complete abelian category).  Under this topology
the points are closed, and when $R$ is left noetherian the
corresponding topological space is noetherian. If $R$ is commutative
(or PI, or FBN) the corresponding topological space is naturally
homeomorphic to the maximal spectrum, equipped with the Zariski
topology. When $R$ is the first Weyl algebra (in characteristic zero)
we obtain a one-dimensional irreducible noetherian topological
space. Comparisons with topologies induced from those on A. L.
Rosenberg's spectra are briefly noted. \endabstract

\address Department of Mathematics, Temple University, Philadelphia,
PA 19122 \endaddress

\email letzter\@math.temple.edu\endemail

\thanks The author's research was supported in part by NSF grants
DMS-9970413 and DMS-0196236. \endthanks

\endtopmatter

\document

\head 1. Introduction \endhead 

One of the fundamental ideas in noncommutative algebraic geometry (see
\cite{\StaVdB} for a recent survey) is that to each noncommutative
ring $R$ there corresponds a ``noncommutative affine space.''
Ideally, such a space should closely reflect the representation theory
of $R$ and should follow a construction mimicking the classical
commutative case. This note, then, is concerned with the
``noncommutative affine space of irreducible representations of $R$.''

Now let {\bf $R$-space\/} denote the set of isomorphism classes of simple
$R$-modules. In \S 2 we equip $R$-space with the {\bf $R$-topology}. It
follows immediately from the definition that points (i.e., simple $R$-modules)
are closed in the $R$-topology, and in \S 3 we prove that the $R$-topology is
noetherian if $R$ is a left noetherian ring. Also, if $R$ is commutative (or
PI, or FBN) then $R$-space, equipped with the $R$-topology is naturally
homeomorphic to $\max R$, equipped with the Zariski topology; see \S 4.

When $R$ is the first Weyl algebra (in characteristic zero), $R$-space
is a one-dimensional, irreducible, noetherian topological space; see
(4.3--4). In particular, the $R$-topology can distinguish between Weyl
algebras and simple Artinian rings (whose corresponding spaces are
singletons).

Following \cite{\Ros} and \cite{\VdB}, we actually work in a somewhat
more general setting. Let $A$ be a complete abelian category, and let
$A$-space denote the collection of isomorphism classes of simple
objects in $A$. Assume further that $A$-space is a set. In \S 2 we
define the $A$-topology on $A$-space, and in \S 3 we prove that the
$A$-topology is noetherian if $A$ has a noetherian generator. As
before, the points in $A$-space are closed. In \S 5 we provide a brief
comparison of the $A$-topology with the topologies developed by
A. L. Rosenberg in \cite{\Ros}. Our notion of a closed set in
$A$-space is also related to ideas found (e.g.) in \cite{\Gab, \Kra,
\KraSao, \Smi}; see (2.6).

My thanks to Ken Goodearl and Paul Smith for their helpful comments on earlier
drafts of this note. I am also grateful to Paul for making his unpublished
lecture notes on noncommutative algebraic geometry available. The questions
considered here were largely inspired by the lectures and workshops on
noncommutative algebraic geometry I attended at MSRI (in February 2000) during
the special year in noncommutative algebra. Finally, I would like to thank the
referee for several helpful suggestions on clarifying the exposition.

\head 2. $A$-space \endhead

\subhead 2.0 Preliminaries \endsubhead The following notation and assumptions
will remain in effect throughout this note. The reader is referred (e.g.) to
\cite{\GooWar, \McCRobii, \Ste} for basic background information on rings and
categories.

(i) We will use $R$ to denote an associative ring with identity.  We will only
use ``$R$-module'' to mean ``left $R$-module,'' and the category of
$R$-modules will be designated $\Mod R$.

(ii) We will use $A$ to denote a complete abelian category. Recall that an
abelian category is complete if and only if it is closed under products
\cite{\Ste, IV.8.3}. Our primary motivating examples of complete abelian
categories are Grothendieck categories \cite{\Ste, X.4.4} and $\Mod R$.

(iii) Following \cite{\VdB}, we will refer to the objects in $A$ as {\bf
$A$-modules\/}, and we will employ the terminology of modules when
appropriate. 

(iv) Let $M$ be an $A$-module, and let $S$ be a set of
submodules of $M$. Set
$$
\bigcap_{N \in S}N \; = \; \text{kernel of the natural map} \; M \;
\rightarrow \; \prod_{N \in S} M/N .
$$
Using this definition, the proof of the Schreier Refinement Theorem (cf, e.g.,
\cite{\GooWar, 3.10}) can be readily adapted to series of $A$-modules.

(v) Let {\bf $A$-space\/} denote the collection of isomorphism classes of
simple $A$-modules (i.e., simple objects in $A$). We will assume for the
remainder of this note that $A$-space is a set. We will use {\bf $R$-space\/}
to denote $(\Mod R)$-space.

(vi) For each $p \in \Aspace$, let $N_p$ denote a chosen representative
$A$-module in $p$. (The topological structure of $A$-space described below
will not depend on these choices.) We will use $[N]$ to denote the isomorphism
class in $A$ of an $A$-module $N$.

\subhead 2.1 \endsubhead Define a subset $X$ of $A$-space to be an {\bf
algebraic set\/} if the isomorphism class of each simple subquotient of
$$
\prod_{p \in X}N_p
$$
is contained in $X$.

\remark{{\bf 2.2} Remark} Suppose that $R$ is commutative. We will see in
(4.1) that the algebraic sets in $R$-space correspond exactly to the Zariski
closed subsets of $\max R$. \endremark

\proclaim{2.3 Theorem} $A$-space, with the closed sets defined to be
the algebraic sets, is a topological space. \endproclaim

\demo{Proof} It is immediately evident that $\emptyset$ and $A$-space
are algebraic sets. Also, it is easy to verify that the intersection of
an arbitary collection of algebraic sets is an algebraic set. Now
suppose that $X_1$ and $X_2$ are algebraic sets, and let $X = X_1 \cup
X_2$. Then
$$
\prod _{p \in X} N_p \; \cong \; \left(\prod_{p \in X_1}N_p\right) \times
\left(\prod_{p \in X_2 \setminus X_1}N_p\right).
$$
It now follows from the Schreier Refinement Theorem that $X$ is an
algebraic set. \qed \enddemo

\subhead 2.4 \endsubhead For convenience, we will refer to the topology
defined in (2.3) as the {\bf $A$-topology}.

\subhead 2.5 \endsubhead For every (two-sided) ideal $I$ of $R$, let
$$
v(I)  \; = \;  \{ p \in \Rspace  \; \colon \; I \subseteq \ann N_p \}.
$$

(i) Generalizing the well-known terminology for primitive ideals,
define the {\bf Jacobson topology\/} on the set $R$-space to be the
topology in which the $v(I)$ are the closed sets. It is easy to see
that the $R$-topology is a refinement of the Jacobson topology.  When
$R$ is commutative, we will refer to the Jacobson topology as the {\bf
Zariski topology\/}.

(ii) Let $\prim R$ denote the (left) primitive spectrum of $R$,
equipped with the usual Jacobson topology. When both $R$-space and
$\prim R$ are equipped with the Jacobson topology, the map
$$
\pi \, \colon \, \Rspace \; @> \; p \mapsto \ann N_p \; >> \; \prim R 
$$
is a closed and continuous surjection.  With respect to the
$R$-topology on $R$-space and the Jacobson topology on $\prim R$, it
follows from (i) that $\pi$ is continuous. When $R$ is commutative,
$\pi$ is a continuous bijection from $R$-space onto $\max R$ (equipped
with the classical Zariski topology); bicontinuity in this case will
follow from (4.1).

\remark{{\bf 2.6} Remarks} (i) In \cite{\Gab}, a full subcategory of an
abelian category is called ``closed'' when it is closed under direct limits
and subquotients; in \cite{\Smi} such subcategories are termed ``weakly
closed.''

(ii) The notion of algebraic set presented in (2.1) is also similar to ideas
developed in \cite{\Kra, \KraSao}.

\endremark

\head 3. Noetherianity \endhead 

Recall the notation of (2.0). Assume in this section that $A$-space is
equipped with the $A$-topology and that $R$-space is equipped with the
$R$-topology.

\subhead 3.1 \endsubhead (i) For each $A$-module $M$, let $S(M)$ denote 
the set of isomorphism classes of simple subquotients of $M$, and let $V(M)$
denote the closure of $S(M)$ in $A$-space.

(ii) For each closed subset $X$ of $A$-space, let
$$
M(X)  \; = \;  \prod_{p \in X}N_p.
$$
Of course, $M(X)$ is determined up to isomorphism only by the closed
set $X$ and not by the choices of $N_p$ for $p \in X$. It follows
from the definitions, since $X$ is closed, that $S(M(X)) = V(M(X)) =
X$.

(iii) If $M$ is an $A$-module we will refer to $M(V(M))$ as the {\bf
  radical\/} of $M$, denoted $\sqrt{M}$. It again follows from the definitions
that
$$
\sqrt{\sqrt{M}} \; \cong \; \sqrt{M}  \quad \text{and} \quad V(M) \; = \;
V\left(\sqrt{M}\right).
$$

\proclaim{3.2 Lemma} Let $X$ be a closed subset of $A$-space, and let $M =
M(X)$. Assume there exists an $A$-module $E$ with the following property: For
each $p \in X$ there is an epimorphism $g_p\colon E \rightarrow N_p$.  Let
$g\colon E \rightarrow M(X)$ denote the product morphism, and let $F$ denote
the image of $g$ in $M(X)$. Then $S(F) = V(F) = X$, and $\sqrt{F} \cong
M$. \endproclaim

\demo{Proof} First, $S(F) \subseteq V(F) \subseteq V(M) = X$, since $F$ is a
submodule of $M$. On the other hand, if $p \in X$ then $N_p$ is isomorphic to
a subquotient of $F$. Hence $X \subseteq S(F)$. Therefore, $X = S(F)$, and
$\sqrt{F} \cong M$. \qed \enddemo

\proclaim{3.3 Theorem} Assume there exists a noetherian $A$-module $E$ that
maps epimorphically onto each simple $A$-module. Then $A$-space is a
noetherian topological space.  \endproclaim

\demo{Proof} Let
$$
X \; = \; X_1 \; \supseteq \; X_2 \; \supseteq \; X_3 \; \supseteq \;
\cdots
$$
be a descending chain of closed subsets of $A$-space. 

For every $p \in X$, choose an epimorphism $g_p \colon E \rightarrow
N_p$. For all $i = 1,2,3,\ldots$, let
$$
E  @>  g_i  >> \prod _{p \in X_i}N_p  \; = \;
M(X_i)
$$
be the product map, and let $M_i$ denote the image of $g_i$ in
$M(X_i)$. Observe that there is an epimorphism
$$
M_i \; \rightarrow \; M_{i+1},
$$
for each $i = 1,2,\ldots$. However, since $M_1$ is noetherian, there
exists some positive integer $t$ for which
$$
M_t \; \cong \; M_{t+1} \; \cong \; \cdots .
$$
Therefore, by (3.2),
$$
X_t \; = \; X_{t+1} \; = \; \cdots .
$$
The theorem follows. \qed \enddemo

\proclaim{3.4 Corollary} If $A$ possesses a noetherian generator then
$A$-space is noetherian, and if $R$ is left noetherian then $R$-space
is noetherian. \qed\endproclaim

\subhead 3.5 \endsubhead Assume that $A$-space is noetherian.
Recalling the standard elementary notions of algebraic geometry, we
see that every algebraic subset is a finite union of irreducible
components. Moreover, we can define the dimension of an algebraic
subset to be the supremum of the lengths of the chains of its
irreducible subsets.

\subhead 3.6 \endsubhead (My thanks to Paul Smith for the following remark;
cf\. \cite{\SmiZan, \S 7}.) We can see as follows that $V(M)$ may be strictly
larger than $S(M)$. Let $k$ be a field, let $\lambda \in k$ be a nonzero
nonroot of unity, and assume that $R = k\{x,y\}/\langle xy - \lambda yx
\rangle$. We can regard the commutative polynomial ring $k[x]$ either as a
subalgebra of $R$ or as an $R$-module on which $y$ acts trivially and $x$ acts
by left multiplication. For $\mu \in k$, set $K(\mu) = R/\langle y, x - \mu
\rangle$, viewed as a $1$-dimensional simple left $R$-module on which $y$ acts
trivially and $x$ acts as multiplication by $\mu$. Set
$$
M \; = \; R\otimes_{k[x]} K(1) \; \cong \; R/R.(x-1) .
$$
It is not hard to see that
$$
S(M) \; = \; \left\{\left[K(\lambda^i)\right] \; \colon \; i = 0,1,2,\ldots
\right\} .
$$
Now let $P = \sqrt{M}$. It is not hard to verify that $P$ contains an
isomorphic copy of the left $R$-module $k[x]$, and it follows, for all
$\mu \in k$, that $K(\mu)$ is isomorphic to a simple $R$-module
quotient of $P$. Therefore, $V(M)$ strictly contains $S(M)$, and
$S(M)$ is not closed in the $R$-topology on $R$-space.

\head 4. When is the $R$-topology equivalent to the Jacobson topology?
\endhead

If the $R$-topology and Jacobson topology coincide on $R$-space, then
every primitive factor of $R$ must have exactly one simple faithful
module (up to isomorphism). In this section we first consider partial
converses to this conclusion, for three specific classes of rings
whose primitive factors are simple artinian: Commutative rings, PI
rings, and FBN rings. We then show that the $R$-topology is a strict
refinement of the Jacobson topology when $R$ is the first Weyl algebra
over a field of characteristic zero.

We retain the notation of the preceding sections.

\proclaim{4.1 Proposition} Suppose that $R$ is a PI ring. Then the
$R$-topology and Jacobson topology coincide on $R$-space. In
particular, when $R$ is commutative the Zariski topology and
$R$-topology coincide on $R$-space. \endproclaim

\demo{Proof} Let $X$ be a subset of $R$-space closed under the
$R$-topology, and let
$$
I \; = \; \ann M(X) \; = \; \bigcap_{p \in X} \ann N_p .
$$
By (2.5i), it suffices to prove that $X = v(I)$. We will assume,
without loss of generality, that $I = 0$, and we will prove that $X$
is equal to all of $R$-space. By Kaplansky's Theorem, there exists a
positive integer $n$ such that $R/\ann N$ is isomorphic to a submodule
of $\oplus_{i=1}^n N$, for all simple $R$-modules $N$. Therefore, $R$
is isomorphic to a submodule of
$$
\prod_{p \in X} \left( \bigoplus_{i=1}^n N_p \right) \; \cong \;
\bigoplus_{i=1}^n M(X) ,
$$
and so every simple $R$-module is isomorphic to a subquotient of
$M(X)$. Thus $R$-space is equal to $X$. \qed\enddemo

\proclaim{4.2 Proposition} Suppose that $R$ is left fully bounded left
noetherian. Then the $R$-topology and Jacobson topology coincide on
$R$-space.  \endproclaim

\demo{Proof} Let $X$ be a subset of $R$-space closed under the
$R$-topology. By (3.2), with $E = {_RR}$, there exists a cyclic
$R$-module $M$ such that $X = S(M)$. By (2.5i), it suffices to prove
that $X = v(\ann M)$, and we will assume without loss of generality
that $\ann M = 0$. Again we must show that $X$ is equal to all of
$R$-space. However, by Cauchon's Theorem (see, e.g., \cite{\GooWar,
8.9}), $R$ embeds as an $R$-module into a finite direct sum of copies
of $M$. Hence every isomorphism class of simple $R$-modules is
contained in $S(M)$, and the proposition follows. \qed\enddemo

We now turn to an example where the Jacobson topology and $R$-topolgy
are distinct. 

\proclaim{4.3 Lemma} Assume that $R$ is a domain with left
Krull dimension equal to $1$. Further suppose that $R$ has infinitely
many pairwise non-isomorphic simple modules. Then $R$-space is a
one-dimensional irreducible topological space. \endproclaim

\demo{Proof} Let $S$ be any infinite collection of maximal left ideals of $R$
for which the simple modules $R/L$, for $L \in S$, are pairwise
non-isomorphic. Because every proper $R$-module factor of $R$ has finite
length, it follows that $\bigcap _{L \in S} L = 0$. Therefore, $R$ embeds as
an $R$-module into every direct product of infinitely many pairwise
non-isomorphic simple $R$-modules.  Consequently, $R$-space itself is the only
infinite closed subset of $R$-space, and so $R$-space is $1$-dimensional and
irreducible.  \qed\enddemo

\subhead 4.4 Example \endsubhead Assume that $k$ is a field of
characteristic zero, and let $R$ denote the first Weyl algebra,
$k\{x,y\}/\langle xy - yx - 1 \rangle$. It is well-known that $R$ is a
simple noetherian domain of left Krull dimension $1$; see, for
example, \cite{\McCRobii}. Moreover, $R$-space is infinite (cf., e.g.,
\cite{\McCRobi}). It now follows from (3.4) and (4.3) that $R$-space
is irreducible, $1$-dimensional, and noetherian. Under the Jacobson
topology, the only closed subsets of $R$-space are $\emptyset$ and
$R$-space itself.

\head 5. Comparisons with topologies relative to A. L. Rosenberg's
spectra \endhead

In this section we assume that the reader is somewhat familiar with the
terminology and notation in \cite{\Ros}, which we will adopt. We will also
continue to use the notation and conventions established in the preceding
sections of this note.

\subhead 5.1 \endsubhead Following \cite{\Ros, \S III.1.2}, $A$-space
can be identified with a subset of $\Spec A$ (as defined in
\cite{\Ros, \S III.1.2}). Several topologies on $\Spec A$ are
considered in \cite{\Ros}, and we can therefore compare the
$A$-topology to their induced (or, relative) topologies on $A$-space.

\subhead 5.2 \endsubhead The topology $\tau$ \cite{\Ros, III.5.1}
induces the discrete topology on $A$-space.

\subhead 5.3 \endsubhead The Zariski topology on $\Spec (\Mod R)$
\cite{\Ros, III.6.3} induces the Jacobson topology on $R$-space. The
central topology \cite{\Ros, III.7.1} on $\Spec (\Mod R)$ is weaker
than this Zariski topology.

\subhead 5.4 \endsubhead (i) In \cite{\Ros, III.7.2}, the topology
$\tau^*$ on $\Spec A$ is defined by declaring the set of supports of
finite type objects in $A$ to be a base.  (The support, in
$\Spec A$, of an object in $A$ is defined in \cite{\Ros, III.5.2}.)
When $A$ is the category of modules over a commutative ring, $\tau^*$
is exactly the Zariski topology on the classical prime spectrum.

(ii) Let $M$ be a finite type object in $A$. Then the set $S(M)$ of ismorphism
classes of simple subquotients of $M$ is closed under the induced $\tau^*$
topology on $A$-space. However, we saw in (3.6) that $S(M)$ need not be closed
under the $A$-topology. Therefore, the $A$-topology and the induced
$\tau^*$-topology are distinct.

(iii) Retain the notation of (3.6). For each non-negative integer $n$,
set 
$$
M(n) \; = \; R\otimes_{k[x]} K(\lambda^n) \; \cong \; R/R.(x-\lambda^n) ,
$$
and let 
$$
S_n \; = \; S(M(n)) \; = \; \{ [K(\lambda^i)] \; \colon \; i = n,n+1,\ldots \}
.
$$
As noted in (ii), each $S_n$ is closed under the induced
$\tau^*$-topology on $R$-space. However, 
$$
S_0 \; \supsetneq \; S_1 \; \supsetneq \; S_2 \; \supsetneq \; \cdots .
$$
Hence, $\tau^*$ need not be a noetherian topology on $R$-space when $R$
is a noetherian ring.

(iv) Suppose that $A$ is a Grothendieck category with a generator $E$
of finite type. Let $X$ be a subset of $A$-space closed under the
$A$-topology. By (3.2), $X = S(F)$ for some quotient $F$ of $E$, and
$F$ is of finite type. Therefore, $X$ is closed under the induced
$\tau^*$ topology. It follows, in this case, that the $\tau^*$
topology is a refinement of the $A$-topology.

\subhead 5.5 \endsubhead The topology $\tau_S$ \cite{\Ros, III.7.3} on
$\Spec A$ is defined to be the weakest topology in which the closure
of a point in $\Spec A$ is equal to its set of specializations.
Therefore, $\tau_S$ reduces to the Zariski topology in the commutative
case, and the points in $A$-space are closed under the induced
$\tau_S$-topology. We do not know whether or not the topology induced
by $\tau_S$ on $A$-space coincides, in general, with the
$A$-topology. We also do not know whether the topology on $A$-space
induced by $\tau_S$ is noetherian when $A$ has a noetherian
generator. A closely related question: Let $M$ be an object in $A$,
and further suppose that $M \in \Spec A$. Must $V(M) = S(M)$?

\enddocument